\documentclass[12pt]{article}

\usepackage{latexsym, amssymb, amsmath, amscd, amsfonts, epsfig, graphicx, colordvi,verbatim,ifpdf}
\usepackage{amsfonts, amsmath, amssymb}
\usepackage{amssymb,amsfonts,amsmath,latexsym,epsfig,cite, psfrag,eepic,color}
\usepackage{amscd,graphics}
\usepackage{latexsym, amssymb,  amsmath,amscd, amsfonts, epsfig, graphicx, colordvi,amsthm}

\usepackage{graphicx}
\usepackage{color}
\usepackage{ifpdf}
\usepackage{fancybox}

\newtheorem{thm}{Theorem}[section]

\newtheorem{rem}[thm]{\it Remarks}
\newtheorem{defi}[thm]{Definition}
\newtheorem{lem}[thm]{Lemma}

\def\qed{\nopagebreak\hfill{\rule{4pt}{7pt}}
\medbreak}

\numberwithin{equation}{section}

\def\qed{\nopagebreak\hfill{\rule{4pt}{7pt}}
\medbreak}

\newlength{\boxedparwidth}
  {\begin{center} \begin{tabular}{|@{\hspace{.315in}}c@{\hspace{.15in}}|}
                  \hline \\ \begin{minipage}[t]{\boxedparwidth}
                  \setlength{\parindent}{.25in}}%
  {\end{minipage} \\ \\ \hline \end{tabular} \end{center}}

\begin{document}
\begin{center}

{ \large\bf  Generalizations of Dyson's Rank on Overpartitions }
\end{center}

\vskip 5mm

\begin{center}
  { Alice X.H. Zhao} \vskip 2mm

College of Science \\[2pt]
Tianjin University of Technology, Tianjin 300384, P.R. China\\[5pt]
 \vskip 2mm
  zhaoxh@email.tjut.edu.cn

\end{center}

\vskip 6mm \noindent {\bf Abstract.} We introduce a statistic on overpartitions called the $\overline{k}$-rank. When there are no overlined parts, this coincides with the $k$-rank of a partition introduced by Garvan.  Moreover, it reduces to the D-rank of an overpartition when $k=2$. The generating function for the $\overline{k}$-rank of overpartitions is  given.
We also establish a relation between the
generating function of  self-3-conjugate overpartitions and the tenth order mock theta functions $X(q)$ and $\chi(q)$.



\section{Introduction}

Dyson’s rank of a partition is defined to be the largest part minus the number of parts \cite{Dyson-1944}. In 1944, Dyson
 conjectured that this partition statistic provided combinatorial interpretations of Ramanujan's congruences
$p(5n+4)\equiv0\pmod5$ and $p(7n+5)\equiv0\pmod7$,
 where $p(n)$ is the number of partitions of $n$.
 Let $N(m,n)$ denote the number of partitions of $n$
with rank $m$. He  also found the following generating function of $N(m,n)$ \cite[eq. (22)]{Dyson-1944}:
\begin{equation}\label{rank}
\sum_{n=0}^\infty N(m,n)q^n=\frac{1}{(q;q)_\infty}\sum_{n=1}^\infty
(-1)^{n-1}q^{n(3n-1)/2+|m|n}(1-q^n).
\end{equation}

Here and in the sequel, we use the standard notation of $q$-series:
\begin{align*}
(a;q)_\infty&:=\prod_{n=0}^\infty(1-aq^n),\
(a;q)_n:=\frac{(a;q)_\infty}{(aq^n;q)_\infty},\\
j(z;q)&:=(z;q)_\infty(q/z;q)_\infty(q;q)_\infty.
\end{align*}
Subsequently, by \eqref{rank} Atkin and Swinnerton-Dyer \cite{Atkin-Swinnerton-Dyer-1954} proved Dyson's conjecture.

Nevertheless, Dyson's rank fails to explain Ramanujan's third
congruence $p(11n+6)\equiv0\pmod{11}$ combinatorially. Thus, Dyson further conjectured that there exists another partition statistic, the crank of a partition, which would likewise give a combinatorial refinement
 of this congruence.

  The crank hypothesized by Dyson was discovered by Andrews and Garvan \cite{Andrews-Garvan-1988}, which is defined as the  the largest part if the partition
 contains no ones, and otherwise as the number of parts larger than the number of
 ones minus the number of ones. Let
 $M(m,n)$ denote  the number of partitions of $n$ with crank $m$. The generating function of $M(m,n)$ was given by \cite[eq. (7.17)]{Garvan-1988},
\begin{equation*}
  \sum_{n=0}^\infty M(m,n)q^n=\frac{1}{(q;q)_\infty}\sum_{n=1}^\infty
(-1)^{n-1}q^{n(n-1)/2+|m|n}(1-q^n).
\end{equation*}

For any positive integer $k$, Garvan \cite{Garvan-1994}
define $N_k(m,n)$ by
\begin{equation}\label{parkg}
\sum_{n=0}^\infty N_k(m,n)q^n=\frac{1}{(q;q)_\infty}\sum_{n=1}^\infty
(-1)^{n-1}q^{n((2k-1)n-1)/2+|m|n}(1-q^n).
\end{equation}
It is easy to see that 
$N_1(m,n)=M(m,n)$ and
 $N_2(m,n)=N(m,n)$.
To give a combinatorial interpretation of $N_k(m,n)$, he introduced the  $k$-rank of a partition by decomposing the Ferrers graph of the partition into successive Durfee squares \cite{Andrews-1979}.
Recall that the  (first) Durfee square of a partition is the largest square within its Ferrers graph and that the $i$-th Durfee square is the Durfee square of
the partition which lies under the ($i-1$)-th Durfee square. For a partition $\pi$ and $j\geq 1$,
  define $n_j(\pi)$ to be the size of the $j$-th Durfee square of $\pi$. If the number of the successive Durfee squares of $\pi$ is less than $l$, then $n_l(\pi)=0$. For convenience, we  define $n_0(\pi)=\infty$.
 The $k$-rank of the partition $\pi$ is defined as the number of columns in the Ferrers graph of $\pi$ which lie to the right of the first Durfee square with length $\leq n_{k-1}(\pi)$ minus the number of parts of $\pi$ that lie below the $(k-1)$-th Durfee square.
 When $k=2$, the $k$-rank coincides with Dyson's rank.
Garvan found that
\begin{thm}{\rm\cite[Theorem 1.12]{Garvan-1994}}\label{krthm}
For $k\geq 2$, $N_k(m,n)$
is the number of partitions of $n$ into at least $k-1$ successive Durfee squares with $k$-rank equal to $m$.
\end{thm}

In the past few years, overpartitions, as a new
component of the theory of partitions, have been heavily studied.  Meanwhile,  many notions  and theorems on ordinary partitions have been extended to overpartitions.

Recall that an overpartition \cite{Corteel-Lovejoy-2004} is a partition in which the first occurrence of each distinct number may be overlined.
The $D$-rank of an overpartition is the largest part  minus the number of parts. Let $\overline{N}(m,n)$ be the number of overpartitions of $n$ with $D$-rank $m$. Lovejoy \cite[Proposition  3.2]{Lovejoy-2005} established the following generating function of $\overline{N}(m,n)$:
\begin{equation}\label{op2}
  \sum_{n=1}^\infty\overline{N}(m,n)q^n=2\frac{(-q;q)_\infty}{(q;q)_\infty}
\sum_{n=1}^\infty\frac{(-1)^{n-1}q^{n^2+|m|n}(1-q^n)}{1+q^n}.
\end{equation}

For any positive integer $k\geq2$, we define $\overline{N}_k(m,n,j)$ by
\begin{equation}\label{k-gen}
   \sum_{n=0}^\infty\sum_{j=0}^\infty\overline{N}_k(m,n,j)a^jq^n
=\frac{(-aq;q)_\infty}{(q;q)_\infty}
\sum_{n=1}^\infty(-1)^{n-1}a^nq^{(k-1)n^2+|m|n}(1-q^n)\frac{(-1/a;q)_n}
{(-aq;q)_n}.
\end{equation}
Setting  $a=0$ and  $a=1$ in \eqref{k-gen} respectively, we can easily deduce 
the following two identities from  \eqref{parkg} and \eqref{op2}:
\begin{align}\label{par-kr}
\overline{N}_k(m,n,0)
=N_k(m,n),\\ \label{op2-1}
\sum_{j=0}^\infty\overline{N}_2(m,n,j)=\overline{N}(m,n).
\end{align}
This inspires us to seek a  combinatorial interpretation of $\overline{N}_k(m,n,j)$. To this end, we first establish the following theorem.
\begin{thm}\label{o-gen}
For  $k\geq 2$, we have
\begin{align}\nonumber
&\sum_{m=-\infty}^\infty\sum_{n=0}^\infty\sum_{j=0}^\infty\overline{N}_k(m,n,j)a^jq^nz^m\\ \label{com-o}
&=\sum_{n_1\geq\cdots\geq n_{k-1}\geq1}
\frac{a^{n_1}q^{{n_1+1\choose 2}+n_2^2+\cdots+n_{k-1}^2}(-1/a;q)_{n_1}}
{(q;q)_{n_1-n_2}\cdots(q;q)_{n_{k-2}-n_{k-1}}
	(zq;q)_{n_{k-1}}(z^{-1}q;q)_{n_{k-1}}}.
\end{align}
\end{thm}

For an overpartition $\lambda$, let $\lambda_1$ denote the largest part of $\lambda$, $l(\lambda)$ denote the number of parts of $\lambda$, and $o(\lambda)$ denote the number of overlined parts of $\lambda$. To interpret \eqref{com-o}
combinatorially, we investigate the vector partitions $(\gamma,\delta,\alpha,\beta)$, where $\gamma$ is a partition of the form
$(l(\gamma),l(\gamma)-1,\ldots,1)$,
 $\delta$ is a partition into distinct non-negative parts less than $l(\gamma)$,
and $\alpha$ and $\beta$ are  partitions with parts at most $l(\gamma)$. We say 
$(\gamma,\delta,\alpha,\beta)$ is a vector partition of $n$ if
 $|\lambda|+|\mu|+|\alpha|+|\beta|=n$.
Let $\mathbb{V}(n)$ denote the set of such vector partitions of $n$. The following lemma can be deduced from the proof of Proposition $2.2$ of \cite{Corteel-Mallet-2007}.
\begin{lem}\label{rel-o-v}
	There exists a bijection  between the vector partitions   $(\gamma,\delta,\alpha,\beta)$ of $n$ counted by $\mathbb{V}(n)$  and overpartitions $\lambda$  of  $n$ such that  $\lambda_1=l(\gamma)+l(\alpha)$, $l(\lambda)=l(\gamma)+l(\beta)$, and  $o(\lambda)=l(\gamma)-l(\delta)$.
\end{lem}	

	Let $\nu(\lambda)$ denote the vector partition corresponding to the overpartition $\lambda$ in the above bijection. We define the $\overline{k}$-rank of an overpartition  as follows.
\begin{defi}
	For an overpartition $\lambda$, let  $\nu(\lambda)=(\gamma,\delta,\alpha,\beta)$. For $k\geq 2$, let $n_{k}(\beta)$ denote the size of the $k$-th successive Durfee square of $\beta$, let $s_k(\alpha)$ denote the	number of parts of $\alpha$ not exceeding $n_k(\beta)$, 
	and let $t_k(\beta)$ denote the number of parts
	of $\beta$ below its $k$-th Durfee square. The
	$\overline{k}$-rank of $\lambda$ is defined  to be $s_{k-2}(\alpha)-t_{k-2}(\beta)$.
\end{defi}
For example, let  $\lambda=(13,10,9,\overline{7},6,\overline{4},4,4,3,1,1,1)$. From the proof of Lemma \ref{rel-o-v}, we will see that $\nu(\lambda)=(\gamma,\delta,\alpha,\beta)$, where  $\gamma=(6,5,4,3,2,1)$, $\delta=(5,3,1,0)$, $\alpha=(5,5,4,2,1,1,1)$   and $\beta=(4,4,3,1,1,1)$.
Note that $n_1(\beta)=3,\ n_2(\beta)=1,\ n_3(\beta)=1$. Then we have $s_3(\alpha)=3$
and $t_3(\beta)=1$, thus the $\overline{5}$-rank of $\lambda$ is $s_3(\alpha)-t_3(\beta)=3-1=2$.

It is worth mentioning that the $\overline{k}$-rank coincides with the D-rank when $k=2$.	For an overpartition $\lambda$ with $\nu(\lambda)=(\gamma,\delta,\alpha,\beta)$, since $n_{0}(\beta)=\infty$,
 $s_0(\alpha)=l(\alpha)$ and $t_0(\beta)=l(\beta)$, we have that the  $\overline{2}$-rank of $\lambda$ is $l(\alpha)-l(\beta)$. Thus, from Lemma \ref{rel-o-v}, we can deduce that the $\overline{2}$-rank of $\lambda$ is equal to  $\lambda_1-l(\lambda)$, i.e. the D-rank of $\lambda$.
 
With the aid of the $\overline{k}$-rank of an overpartition, 
 we give a combinatorial interpretation of $\overline{N}_k(m,n,j)$.

\begin{thm}\label{com-rank}
	For  $k\geq 2$, $\overline{N}_{k}(m,n,j)$ is the number of overpartitions $\lambda$ of $n$ into at least $k-2$ successive Durfee squares in $\beta$ with $j$ overlined parts and $\overline{k}$-rank $m$, where $\beta$ is the fourth component of  $\nu(\lambda)$.
\end{thm}

 In view of \eqref{par-kr}, we will see that the $j=0$ case of Theorem \ref{com-rank} reduces to Theorem \ref{krthm}. From \eqref{op2-1} and Theorem \ref{com-rank}, we can also find that
$\overline{N}_{2}(m,n,j)$ is a dissection of $\overline{N}(m,n)$ according to the number of overlined parts.

From \eqref{k-gen},	it is easily seen that
\begin{equation}\label{cj}
\overline{N}_k(m,n,j)=\overline{N}_k(-m,n,j).
\end{equation}
For an overpartition $\lambda$, let 
$\nu(\lambda)=(\gamma,\delta,\alpha,\beta)$. We find that \eqref{cj} can be  proved combinatorially  by interchanging the parts enumerated by $s_{k-2}(\alpha)$ and $t_{k-2}(\beta)$. We call this operation as the $k$-conjugation, which is essentially the same as the $k$-conjugation
introduced by Corteel, Lovejoy and Mallet \cite[p. 1616]{Corteel-Lovejoy-Mallet}. If an overpartition is fixed by $k$-conjugation, then it is called a   self-$k$-conjugate overpartition.  Corteel, Lovejoy and Mallet derived an identity related to the generating function of self-$k$-conjugate
overpartitions \cite[eq (1.8)]{Corteel-Lovejoy-Mallet}:
\begin{align}\nonumber
&\sum_{n_1\geq\ldots \geq n_{k-1}\geq 0}\frac{q^{{n_1+1\choose 2}+n^2_2+\cdots+n^2_{k-1}}a^{n_1}(-1/a;q)_{n_1}}
{(q;q)_{n_1-n_2}\cdots(q;q)_{n_{k-2}-n_{k-1}}(q^2;q^2)_{n_{k-1}}}\\ \label{sk-con-gen}
&\hskip2cm=\frac{(-aq;q)_\infty}{(q;q)_\infty}\sum_{n=-\infty}^\infty\frac{(-1/a;q)_n(-1)^na^nq^{(k-1)n^2+{n+1\choose 2}}}{(-aq;q)_n}.
\end{align}

From the $k=3$ case of \eqref{sk-con-gen},  we  derive the following identity equivalent to \cite[eq. (1.2)]{Lovejoy-Osburn-2013}. Especially, this identity establishes a relation between the
generating function for  self-3-conjugate overpartitions and the tenth order mock theta functions $X(q)$ and $\chi(q)$ \cite{Ramanujan-1988}, where
\begin{align*}
X(q)=\sum_{n=0}^\infty\frac{(-1)^nq^{n^2}}{(-q;q)_{2n}} \text{\ \ \ and\ \ \ }
\chi(q)=\sum_{n=0}^\infty \frac{(-1)^nq^{(n+1)^2}}{(-q;q)_{2n+1}}.
\end{align*}

\begin{thm} \label{eqmock}
	We have
	\begin{equation*}
	\sum_{n_1\geq n_2\geq0}\frac{(-1;q)_{n_1}q^{{n_1+1\choose 2}+n^2_2}}
	{(q;q)_{n_1-n_2}(q^2;q^2)_{n_2}}=\frac{2(-q;q)_\infty X(q) }{(q;q^5)(q^4;q^5)_\infty}-
	\frac{2(-q;q)_\infty \chi(q) }{(q^2;q^5)(q^3;q^5)_\infty}-\frac{(q;q)_\infty}{(-q;q)_\infty}.
	\end{equation*}
\end{thm}

This paper is organized as follows.  In Section $2$,  we give a proof of  Theorem \ref{o-gen} by utilizing Andrews' multiple series transformation.  Section 3 is devoted   to giving a combinatorial interpretation of $\overline{N}_k(m,n)$, which is stated in Theorem \ref{com-rank}.    In Section 4, we aim to show the relation between the generating function of self-3-conjugate overpartitions and  mock theta functions as stated in Theorem \ref{eqmock}.

\section{Proof of Theorem \ref{o-gen}}

 In this section, we prove that Theorem \ref{o-gen} follows from Andrews' $k$-fold generalization of $q$-Whipple's theorem \cite{Andrews-1975}, which is stated as follows.
For any positive integer $k$,
\begin{align}\nonumber
&_{2k+4}\phi_{2k+3}
\setlength\arraycolsep{1mm}\left(\begin{array}{ccccccccc}a,&qa^{\frac{1}{2}},
&-qa^{\frac{1}{2}},&b_1,&c_1,&\ldots,&b_k,&c_k,&q^{-N}\\
&a^{\frac{1}{2}},&-a^{\frac{1}{2}},&aq/b_1,&aq/c_1,&\ldots,
&aq/b_k,&aq/c_k,&aq^{N+1}\end{array};q,\frac{a^k
q^{k+N}}{b_1\cdots b_kc_1\cdots c_k}\right)\\[5pt] \nonumber
&=\frac{(a q;q)_{N}\left(a q / b_{k} c_{k};q\right)_{N}}{\left(a q / b_{k};q\right)_{N}\left(a q / c_{k};q\right)_{N}} \sum_{m_{1}, \ldots, m_{k-1} \geq 0} \frac{\left(a q / b_{1} c_{1};q\right)_{m_{1}}\left(a q / b_{2} c_{2};q\right)_{m_{2}} \cdots\left(a q / b_{k-1} c_{k-1};q\right)_{m_{k-1}}}{(q;q)_{m_{1}}(q;q)_{m_{2}} \cdots(q;q)_{m_{k-1}}}\\ \nonumber
&\ \ \ \times \frac{\left(b_{2};q\right)_{m_{1}}\left(c_{2};q\right)_{m_{1}}
\left(b_{3};q\right)_{m_{1}+m_{2}}\left(c_{3};q\right)_{m_{1}+m_{2}} \cdots\left(b_{k};q\right)_{m_{1}+\cdots+m_{k-1}}}{\left(a q / b_{1};q\right)_{m_{1}}\left(a q / c_{1};q\right)_{m_{1}}\left(a q / b_{2};q\right)_{m_{1}+m_{2}}\left(a q / c_{2};q\right)_{m_{1}+m_{2}} \cdots\left(a q / b_{k-1};q\right)_{m_{1}+\cdots+m_{k-1}}} \\ \nonumber
&\ \ \ \times \frac{\left(c_{k};q\right)_{m_{1}+\cdots+m_{k-1}}}{\left(a q / c_{k-1};q\right)_{m_{1}+\cdots+m_{k-1}}} \cdot \frac{\left(q^{-N};q\right)_{m_{1}+\cdots+m_{k-1}}}{\left(b_{k} c_{k} q^{-N} / a;q\right)_{m_{1}+\cdots+m_{k-1}}} \\ \label{k-fold}
&\  \ \ \times \frac{(a q)^{m_{k-2}+2 m_{k-3}+\cdots+(k-2) m_{1}} q^{m_{1}+\cdots+m_{k-1}}}{\left(b_{2} c_{2}\right)^{m_{1}}\left(b_{3} c_{3}\right)^{m_{1}+m_{2}} \cdots\left(b_{k-1} c_{k-1}\right)^{m_{1}+\cdots+m_{k-2}}},
\end{align}
here
\begin{equation*}
_{r+1}\phi_r\left(\begin{array}{cccc}a_0,&a_1,&\ldots,&a_r\\
&b_1,&\ldots,&b_r\end{array};q,z\right)=\sum_{n=0}^\infty
\frac{(a_0;q)_n(a_1;q)_n\cdots(a_r;q)_nz^n}{(q;q)_n(b_1;q)_n\cdots(b_r;q)_n}.
\end{equation*}
\noindent{\it Proof of Theorem \ref{o-gen}\ } In \eqref{k-fold}, let $a=1, b_1=c_1^{-1}=z$, then set $b_{k}=-1/a$ and let $N$ and all the other $b_i$ and $c_i$ tend to infinity. After simplification, we obtain
\begin{align} \nonumber
   &1+\sum_{n=1}^\infty\frac{(-1)^na^n(-1/a;q)_nq^{(k-1)n^2+n}(1+q^n)
   (1-z)(1-z^{-1})}
{(-aq;q)_n(1-zq^n)(1-z^{-1}q^n)}\\[5pt] \nonumber
&=\frac{(q;q)_\infty}{(-aq;q)_\infty}\sum_{m_1,\ldots,m_{k-1}\geq0}
\frac{q^{{m_1+\cdots+m_{k-1}+1\choose 2}+(m_1+\cdots+m_{k-2})^2+\cdots+m_1^2}}
{(q;q)_{m_2}\cdots(q;q)_{m_{k-1}}}\\ \label{k-fold-1}
&\hskip 5cm\times\frac{a^{m_1+\cdots+m_{k-1}}(-1/a;q)_{m_1+\cdots+m_{k-1}}}{(zq;q)_{m_1}(z^{-1}q;q)_{m_1}}.
\end{align}
Multiplying both sides of \eqref{k-fold-1} by $(-aq;q)_\infty/(q;q)_\infty$, and setting $n_{k-i}=m_1+m_2+\cdots+m_i$ for $1\leq i\leq k-1$, we have
\begin{align}\nonumber
   &\frac{(-aq;q)_\infty}{(q;q)_\infty}
   \left(1+\sum_{n=1}^\infty\frac{(-1)^na^n(-1/a;q)_nq^{(k-1)n^2+n}(1+q^n)
   (1-z)(1-z^{-1})}
{(-aq;q)_n(1-zq^n)(1-z^{-1}q^n)}\right)\\[5pt]\label{k-fold-2}
=&\sum_{n_1\geq\cdots\geq n_{k-1}\geq0}
\frac{a^{n_1}q^{{n_1+1\choose 2}+n_2^2+\cdots+n_{k-1}^2}(-1/a;q)_{n_1}}
{(q;q)_{n_1-n_2}\cdots(q;q)_{n_{k-2}-n_{k-1}}
(zq;q)_{n_{k-1}}(z^{-1}q;q)_{n_{k-1}}}.
\end{align}

Invoking Garvan's Lemma \cite[Lemma 3.6]{Garvan-1994}
 \[\frac{
   q^n(1-z)(1-z^{-1})}{(1-zq^n)(1-z^{-1}q^n)}=
1-\frac{1-q^n}{1+q^{n}}\sum_{m=0}^\infty z^mq^{mn}-
   \frac{1-q^n}{1+q^n}\sum_{m=1}^\infty z^{-m}q^{mn},\]
we find that the left side of \eqref{k-fold-2} equals
\begin{align*}
   &\frac{(-aq;q)_\infty}{(q;q)_\infty}
   \left\{1+\sum_{n=1}^\infty(-1)^na^nq^{(k-1)n^2}(1+q^n)\frac{(-1/a;q)_n}
{(-aq;q)_n}\right.\\[5pt]
&\hskip 5cm \times\left.\left(1-\frac{1-q^n}{1+q^{n}}\sum_{m=0}^\infty z^mq^{mn}-
   \frac{1-q^n}{1+q^n}\sum_{m=1}^\infty z^{-m}q^{mn}\right)\right\}\\
   &=\frac{(-aq;q)_\infty}{(q;q)_\infty}\left(1+
   \sum_{n=1}^\infty(-1)^na^nq^{(k-1)n^2}(1+q^n)\frac{(-1/a;q)_n}
{(-aq;q)_n}\right)\\
&\hskip0.5cm +\frac{(-aq;q)_\infty}{(q;q)_\infty}
\sum_{n=1}^\infty(-1)^{n-1}a^nq^{(k-1)n^2}(1-q^n)\frac{(-1/a;q)_n}
{(-aq;q)_n}\left(\sum_{m=0}^\infty z^mq^{mn}+
  \sum_{m=1}^\infty z^{-m}q^{mn}\right).
\end{align*}

It is easy to see that
\begin{align*}
  \sum_{n=1}^\infty(-1)^na^nq^{(k-1)n^2}\frac{(-1/a;q)_n}
{(-aq;q)_n}
=\sum_{n=-\infty}^{-1}(-1)^na^nq^{(k-1)n^2+n}\frac{(-1/a;q)_n}
{(-aq;q)_n}.
\end{align*}
Thus we have
\begin{align}\nonumber
 &\frac{(-aq;q)_\infty}{(q;q)_\infty}\left(1+
   \sum_{n=1}^\infty(-1)^na^nq^{(k-1)n^2}(1+q^n)\frac{(-1/a;q)_n}
{(-aq;q)_n}\right)\\ \nonumber
&=\frac{(-aq;q)_\infty}{(q;q)_\infty}
\sum_{n=-\infty}^\infty(-1)^na^nq^{(k-1)n^2+n}\frac{(-1/a;q)_n}
{(-aq;q)_n}\\ \label{k-fold-4}
&=\sum_{n_1\geq\cdots\geq n_{k-2}\geq0}
\frac{a^{n_1}q^{{n_1+1\choose 2}+n_2^2+\cdots+n_{k-2}^2}(-1/a;q)_{n_1}}
{(q;q)_{n_1-n_2}\cdots(q;q)_{n_{k-3}-n_{k-2}}
(q;q)_{n_{k-2}}},
\end{align}
where the last equality follows from \cite[eqs. (1.1) and (6.1)]{Corteel-Mallet-2007}
\begin{align*}
  &\frac{(-aq;q)_\infty}{(q;q)_\infty}
\sum_{n=-\infty}^\infty(-1)^na^nq^{kn^2+(k-i+1)n}\frac{(-1/a;q)_n}
{(-aq;q)_n}\\
&=\sum_{n_1\geq\cdots\geq n_{k-1}\geq0}
\frac{a^{n_1}q^{{n_1+1\choose 2}+n_2^2+\cdots+n_{k-1}^2+n_i+\cdots+n_{k-1}}(-1/a;q)_{n_1}}
{(q;q)_{n_1-n_2}\cdots(q;q)_{n_{k-2}-n_{k-1}}
(q;q)_{n_{k-1}}}
\end{align*}
with $i=k$ and $k$ replaced by $k-1$.

We note that the last term in \eqref{k-fold-4} corresponds to the part of the sum on the right side of
\eqref{k-fold-2} with $n_{k-1}=0$.
Subtracting this term from both sides of \eqref{k-fold-2} by using \eqref{k-fold-4}, we deduce that
\begin{align*}
 &\sum_{n_1\geq\cdots\geq n_{k-1}\geq1}
\frac{a^{n_1}q^{{n_1+1\choose 2}+n_2^2+\cdots+n_{k-1}^2}(-1/a;q)_{n_1}}
{(q;q)_{n_1-n_2}\cdots(q;q)_{n_{k-2}-n_{k-1}}
(zq;q)_{n_{k-1}}(z^{-1}q;q)_{n_{k-1}}} \\
   &=\frac{(-aq;q)_\infty}{(q;q)_\infty}
\sum_{n=1}^\infty(-1)^{n-1}a^nq^{(k-1)n^2}(1-q^n)\frac{(-1/a;q)_n}
{(-aq;q)_n}\left(\sum_{m=0}^\infty z^mq^{mn}+
  \sum_{m=1}^\infty z^{-m}q^{mn}\right)\\ 
  &=\sum_{m=-\infty}^\infty\sum_{n=0}^\infty
  \sum_{j=0}^\infty\overline{N}_k(m,n,j)a^jq^nz^m.
\end{align*}

 Thus we complete the proof of Theorem \ref{o-gen}. \qed

\section{Combinatorial interpretation of $\overline{N}_k(m,n,j)$}

In this section we mainly give a proof of Theorem \ref{com-rank}. To achieve this, we shall first show 
 Lemma \ref{rel-o-v} with the aid of the generalized Durfee square \cite{Corteel-Mallet-2007} of an overpartition.
  Recall that the size of the generalized Durfee square of an overpartition is defined to be the largest number $N$ such that the number of overlined parts plus the
 number of non-overlined parts greater than or equal to $N$ is at least $N$. For example, the size of the generalized Durfee square of $\lambda=(13,10,9,\overline{7},6,\overline{4},4,4,3,1,1,1)$ is 6.

We now proceed to give a proof of Lamma \ref{rel-o-v}.

\noindent{\it Proof of Lemma \ref{rel-o-v} } 
For a vector partition $(\gamma, \delta,\alpha,\beta)\in \mathbb{V}(n)$, let $N=l(\gamma)$ and let $\alpha'=(\alpha_1',\alpha_2',\ldots,\alpha_N')$ be the conjugate of $\alpha=(\alpha_1,\alpha_2,\ldots,)$. If $\alpha_1$  is less than $N$, we set $\alpha_j'=0$ for $\alpha_1<j\leq N$. We will construct an overpartition $\lambda$ as follows.  We first add  $\gamma$ to $\alpha'$ to get a distinct partition $\sigma=(\alpha_1'+N,\alpha_2'+N-1,\ldots,\alpha_N'+1)$. 
It is easy to see that $l(\sigma)=N$ and $\sigma_i=\alpha_i'+N-i+1\geq N-i+1$ for $1\leq i\leq N$. Next we overline all the parts of $\sigma$. For each part $s$ in $\delta$, we add  $s$ to the part $\sigma_{s+1}$ and remove its overline.   Noting that $\sigma_{s+1}+s\geq N$, after reordering the parts we then get an overpartition $\pi$ with each non-overlined part $\geq N$. Finally, by jointing  $\pi$
and $\beta$ together, we get an overpartition $\lambda$ of $n$ with generalized Durfee square of size $N$.
From the above procedure, we can easily check that the largest part of $\lambda$ is the largest part of $\sigma$, which equals to $l(\alpha)+l(\gamma)$. The number of parts of $\lambda$ is $l(\gamma)+l(\beta)$. And the number of overlined parts of $\lambda$ is $l(\gamma)-l(\delta)$.
 
For example, let $\nu=(\gamma, \delta,\alpha,\beta)\in \mathbb{V}(63)$, where    $\gamma=(6,5,4,3,2,1)$, $\delta=(5,3,1,0)$, $\alpha=(5,5,4,2,1,1,1)$ and
$\beta=(4,4,3,1,1,1)$. Then we have $\alpha'=(7,4,3,3,2,0)$, and
$\sigma=(13,9,7,6,4,1)$. The overpartition $\pi$ follows from $\sigma$ and $\delta$ according to the following procedure:
\begin{align*}
(\sigma,\delta)&=(13+9+7+6+4+1,5+3+1+0)\\
&\Rightarrow (\overline{13}+\overline{9}+\overline{7}+\overline{6}+\overline{4}+\overline{1},5+3+1+0)\\
&\Rightarrow (13+\overline{9}+\overline{7}+\overline{6}+\overline{4}+\overline{1},5+3+1)\\
&\Rightarrow (13+10+\overline{7}+\overline{6}+\overline{4}+\overline{1},5+3)\\
&\Rightarrow (13+10+9+\overline{7}+\overline{4}+\overline{1},5)\\
&\Rightarrow (13+10+9+\overline{7}+6+\overline{4})=\pi.
\end{align*}
Finally, jointing $\pi$ and $\beta$ together, we derive  $\lambda=(13,10,9,\overline{7},6,\overline{4},4,4,3,1,1,1)$, which is an overpartition of 63.

The above process can be uniquely reversed.
For an overpartition $\lambda$ of $n$,  let $j$ be the number of  its overlined parts and $N$ be the size of its generalized Durfee square. 
The vector partition $(\gamma, \delta,\alpha,\beta)$ is constructed  as follows. 

We first decompose $\lambda$ into three partitions $\rho=(\lambda_1,\ldots,\lambda_j)$, $\mu=(\lambda_{j+1},\ldots,\lambda_N)$ and our desired $\beta=(\lambda_{N+1},\ldots,\lambda_{l(\lambda)})$,
where $\lambda_p\ (1\leq p\leq j)$ is the overlined part of $\lambda$ with its overline removed and $\lambda_1>\lambda_2> \cdots > \lambda_j$ and $\lambda_q\ (j+1\leq q\leq l(\lambda))$ is the non-overlined part of $\lambda$ such that
$\lambda_{j+1}\geq\lambda_{j+2}\geq\cdots\geq  \lambda_{l(\lambda)}$. From the definition of the generalized Durfee square of an overpartition, we can check  that $\lambda_N\geq N$ and 
$\lambda_{N+1}\leq N$.
To get the distinct non-negative  partition $\delta$, we initialize $\delta=\emptyset$ and $\sigma=\rho$ and do the following operation to $\lambda_t$ for $t$ from $N$ to $j+1$. 
 Let $\sigma_i$  be the $i$-th  part of $\sigma$.
 For convenience, let $\sigma_0=\infty$.
Let $s$ be the largest number  such that $\lambda_t-t-s+j+1 < \sigma_{s}$. If $s\geq 1$, then we insert the part $\lambda_t-t-s+j+1$ into $\sigma$ right after $\sigma_s$ and add the part $t+s-j-1$ to $\delta$. If $s=0$, then we insert the part  $\lambda_t-t+j+1$ into $\sigma$ before $\sigma_1$. 
 After the operation, we derive a distinct partition $\sigma$  and   a distinct partition $\delta$ with parts lying in $[0,N-1]$. Next we split $\sigma$
into two partitions $\gamma=(N,N-1,\ldots,1)$ and $\eta$, where $\eta_i=\sigma_i+i-1-N$ for $1\leq i\leq N$. Taking the  conjugate of $\eta$, we then get the partition $\alpha$ with  parts $\leq N$. Finally, we get the vector partition  $(\gamma, \delta,\alpha,\beta)\in \mathbb{V}(n)$.

For example, let $\lambda=(13,10,9,\overline{7},6,\overline{4},4,4,3,1,1,1)$ be an overpartition of 63. It is easy to see that $j=2$ and $N=6$. Thus we have $\rho=(7,4)$, $\mu=(13,10,9,6)$ and  $\beta=(4,4,3,1,1,1)$. For each part in $\mu$, we do the following operation
\begin{align*}
(\mu,\rho,\delta)&=(13+10+9+6,7+4,\emptyset)\\
&\Rightarrow(13+10+9,7+4+1,5)\\
&\Rightarrow(13+10,7+6+4+1,5+3)\\
&\Rightarrow(13,9+7+6+4+1,5+3+1)\\
&\Rightarrow(\emptyset,13+9+7+6+4+1,5+3+1+0)=(\emptyset,\sigma,\delta).
\end{align*} 
Then we have $\delta=(5,3,1,0)$, $\gamma=(6,5,4,3,2,1)$, $\eta=(7,4,3,3,2,0)$ and $\alpha=\eta'=(5,5,4,2,1,1,1)$. Hence we get a vector partition  $\nu=(\gamma, \delta,\alpha,\beta)\in \mathbb{V}(63)$.

Therefore we complete the proof of Lemma \ref{rel-o-v}.\qed

We are ready to present the proof of Theorem \ref{com-rank}.

\noindent{\it Proof of Theorem \ref{com-rank}\ } First let us rewrite the sum on the right side  of \eqref{com-o} as
\begin{align}\nonumber
   &\sum_{n_1\geq \cdots \geq n_{k-1}\geq1} q^{{n_1+1\choose 2}}
   \times
   \frac{1}{(1-q^{n_{k-1}+1})\cdots(1-q^{n_1})(1-zq)\cdots(1-zq^{n_{k-1}})}\\ \label{gen-3}
&\hskip0.5cm \times a^{n_1}(-\frac{1}{a};q)_{n_1}\times
q^{n^2_{2}}{n_1 \brack n_2}\times\cdots\times q^{n^2_{k-1}}{n_{k-2}\brack n_{k-1}}
\times\frac{1}{(1-z^{-1}q)\cdots(1-z^{-1}q^{n_{k-1}})}.
\end{align}

The factor $q^{n_1+1\choose2}$ generates the partition $\gamma=(n_1,n_1-1,\ldots,1)$. The term
 \[a^{n_1}(-\frac{1}{a};q)_{n_1}\]
 generates partitions $\delta$ into distinct non-negative parts less than $n_1$, and the power of $a$ counts $n_1-l(\delta)$.

Since
\[{n\brack m}\] is the generating function for partitions into at most $m$ parts each $\leq n-m$, the term
\[\sum_{n_2\geq\cdots\geq n_{k-1}\geq 1}q^{n^2_{2}}{n_1 \brack n_2}\times\cdots\times q^{n^2_{k-1}}{n_{k-2}\brack n_{k-1}}\times
\frac{1}{(1-z^{-1}q)\cdots(1-z^{-1}q^{n_{k-1}})}\]
generates partitions $\beta$ into parts  $\leq n_1$ and at least  $k-2$ successive Durfee squares. The power of $z^{-1}$ keeps track the number of parts that lie below its $(k-2)$-th Durfee square.

The term
\[\frac{1}{(1-q^{n_{k-1}+1})\cdots(1-q^{n_1})(1-zq)\cdots(1-zq^{n_{k-1}})}\]
generates partitions $\alpha$ into parts $\leq n_1$ and the power of $z$ counts the number of parts $\leq n_{k-1}$, where $n_{k-1}$ is the size of the $(k-2)$-th Durfee square of $\beta$.

Thus, \eqref{gen-3} is the generating function of vector partitions $(\gamma,\delta,\alpha,\beta)\in\cup_n \mathbb{V}(n)$.  The power of  $a$ in \eqref{gen-3}
counts $l(\gamma)-l(\delta)$, and the power of $z$ counts the number of parts of $\alpha$ that are $\leq n_{k-1}$  minus the number of parts  of $\beta$ that lie below the $(k-2)$-th Durfee square of $\beta$, where $n_{k-1}$
is the size of $(k-2)$-th Durfee square of $\beta$.
In light of Lemma \ref{rel-o-v}, we can deduce that \eqref{gen-3} also generates  overpartitions $\lambda$  into at least $k-2$ successive Durfee squares in $\beta$, where  $\beta$ is the fourth component of the vector partition $\nu(\lambda)$.
Moreover, the  power of $a$  counts the number of overlined parts of $\lambda$, i.e. $o(\lambda)$, and
  the power of $z$  counts the $\overline{k}$-rank of $\lambda$.

 Therefore Theorem \ref{com-rank} follows from \eqref{com-o} and \eqref{gen-3}.\qed

\begin{rem}	{\rm
For an overpartition $\lambda$, let $\nu(\lambda)=(\gamma,\delta,\alpha,\beta)$. Note that if there are no overlined parts in $\lambda$ (i.e., $\lambda$ is a partition), then from the proof of Lemma \ref{rel-o-v} we can see that $\delta=(l(\gamma)-1, l(\gamma)-2,\ldots,1,0)$ and 
  the $(k-2)$-th Durfee square of  $\beta$ is in fact the  $(k-1)$-th Durfee square of  $\lambda$. And from the definition of $\overline{k}$-rank, it is easy to see that the  $\overline{k}$-rank of $\lambda$ is actually the
$k$-rank of $\lambda$. Recall that $\overline{N}_k(m,n,0)=N_k(m,n)$. Consequently,  the $j=0$ case of Theorem \ref{com-rank} coincides with Theorem \ref{krthm}.}
\end{rem}

\section{Proof of Theorem \ref{eqmock}}

This section is devoted to showing the relation between the generating function of self-3-conjugate overpartitions and  the mock theta functions as stated in Theorem \ref{eqmock}. To this end, we mainly express both of them in terms of
 the Appell-Lerch sum $m(x,q,z)$, where
\[m(x,q,z)=\frac{1}{j(z;q)}\sum_{n=-\infty}^\infty\frac{(-1)^nq^{n\choose 2}z^n}
{1-q^{n-1}xz}.\]

Hickerson and Mortenson \cite{Hickerson-Mortenson-2014} showed that
\begin{align}
\label{tenord-1}
X(q)&=2m(-q^2,q^5,q^4)
-\frac{j(q^3;q^{10})j(q^5;q^{10})}{j(q;q^{5})},\\ \label{tenord-2}
\chi(q)&=2m(-q,q^5,q^2)
+q\frac{j(q;q^{10})j(q^5;q^{10})}{j(q^2;q^{5})}.
\end{align}

Now we proceed to present our proof of Theorem \ref{eqmock}.

\noindent{\it Proof of Theorem \ref{eqmock}\ }
 Setting $a=1$ and $k=3$ in  	\eqref{sk-con-gen}, it follows that
 \allowdisplaybreaks
\begin{align*}
 \sum_{n_1\geq n_2\geq0}\frac{(-1;q)_{n_1}q^{{n_1+1\choose 2}+n^2_2}}
 {(q;q)_{n_1-n_2}(q^2;q^2)_{n_2}}
 =\frac{2(-q;q)_\infty}{(q;q)_\infty}\sum_{n=-\infty}^\infty\frac{(-1)^nq^{5{n\choose 2}+3n}}{1+q^{n}}.
 \end{align*}
 Since
\[\frac{1}{1-x}=\frac{1+x+\cdots+x^{k-1}}{1-x^{k}},\]
we have
\begin{align}\nonumber
 &\sum_{n_1\geq n_2\geq0}\frac{(-1;q)_{n_1}q^{{n_1+1\choose 2}+n^2_2}}
 {(q;q)_{n_1-n_2}(q^2;q^2)_{n_2}}\\ \nonumber
&=\frac{2(-q;q)_\infty}{(q;q)_\infty}\sum_{n=-\infty}^\infty\frac{(-1)^nq^{5{n\choose 2}+3n}(1-q^n+q^{2n}-q^{3n}+q^{4n})}{1+q^{5n}}\\ \nonumber
&=\frac{2(-q;q)_\infty}{(q;q)_\infty}\left(\sum_{n=-\infty}^\infty\frac{(-1)^nq^{5{n\choose 2}+3n}}{1+q^{5n}}-\sum_{n=-\infty}^\infty\frac{(-1)^nq^{5{n\choose 2}+4n}}{1+q^{5n}}+\sum_{n=-\infty}^\infty\frac{(-1)^nq^{5{n+1\choose 2}}}{1+q^{5n}}\right.\\ \nonumber
&\hskip3cm\left.-\sum_{n=-\infty}^\infty\frac{(-1)^nq^{5{n\choose 2}+6n}}{1+q^{5n}}+\sum_{n=-\infty}^\infty\frac{(-1)^nq^{5{n\choose 2}+7n}}{1+q^{5n}}\right)\\ \nonumber
&=\frac{2(-q;q)_\infty}{(q;q)_\infty}\left(2\sum_{n=-\infty}^\infty\frac{(-1)^nq^{5{n\choose 2}+3n}}{1+q^{5n}}-2\sum_{n=-\infty}^\infty\frac{(-1)^nq^{5{n\choose 2}+4n}}{1+q^{5n}}+\sum_{n=-\infty}^\infty\frac{(-1)^nq^{5{n+1\choose 2}}}{1+q^{5n}}\right)\\ \label{2cm-1}
&=\frac{2(-q;q)_\infty}{(q;q)_\infty}\left(2j(q^3;q^5)m(-q^2,q^5,q^3)-
2j(q^4;q^5)m(-q,q^5,q^4)+\frac{(q^5;q^5)_\infty^3}{j(-1;q^5)}\right),
\end{align}
where  the last equality follows from \cite[p. 1]{Ramanujan-1988}
\[\sum_{n=-\infty}^\infty\frac{(-1)^nq^{n+1\choose 2}}{1-zq^n}
=\frac{(q;q)_\infty^3}{j(z;q)}.
\]

Recall that for generic $x,z_0,z_1\in\mathbb{C}^*$,
the following identity \cite[eq. (3.7)]{Hickerson-Mortenson-2014} holds
\begin{equation}\label{mmtrans}
m(x,q,z_1)-m(x,q,z_0)=\frac{z_0(q;q)_\infty^3j(z_1/z_0;q)j(xz_0z_1;q)}
{j(z_0;q)j(z_1;q)j(xz_0;q)j(xz_1;q)}.
\end{equation}
From \eqref{mmtrans},  we can deduce that
\begin{align}\label{mdif-1}
 m(-q^2,q^5,q^3)-m(-q^2,q^5,q^4)&=\frac{-(q^5;q^5)_\infty^3}{j(q^3;q^5)j(-q^5;q^5)},\\ \label{mdif-2}
 m(-q,q^5,q^4)-m(-q,q^5,q^2)&=\frac{(q^5;q^5)_\infty^3}{j(q^4;q^5)j(-q^5;q^5)}.
\end{align}

Substituting \eqref{tenord-1}, \eqref{tenord-2},
  \eqref{mdif-1} and \eqref{mdif-2} into \eqref{2cm-1} yields
\begin{align*}
  &\sum_{n_1\geq n_2\geq0}\frac{(-1;q)_{n_1}q^{{n_1+1\choose 2}+n^2_2}}
  {(q;q)_{n_1-n_2}(q^2;q^2)_{n_2}}\\
&=\frac{2(-q;q)_\infty j(q^3;q^5)}{(q;q)_\infty}X(q)-
\frac{2(-q;q)_\infty j(q^4;q^5)}{(q;q)_\infty}\chi(q)+\frac{2(-q;q)_\infty}{(q;q)_\infty}\\
&\quad\times\left(\frac{j(q^3;q^5)j(q^3;q^{10})
j(q^5;q^{10})}{j(q;q^{5})}
+q\frac{j(q^4;q^5)j(q;q^{10})j(q^5;q^{10})}{j(q^2;q^{5})}
-\frac{3(q^5;q^5)_\infty^3}{j(-1;q^5)}
\right).
\end{align*}
To finish the proof of Theorem \ref{eqmock}, we just need to show that the series in the  brackets on  the right hand side is equal to $-(q;q)^2_\infty/2(-q;q)^2_\infty$. Note that this is  an identity between modular forms, thus we can easily verify it with a finite computation by using Garvan's MAPLE program, thetaids \cite{Frye-Garvan-2019}, available at
https://qseries.org/fgarvan/qmaple/thetaids.
Hence we complete the proof of Theorem \ref{eqmock}.
\qed

    \vspace{0.5cm}
\noindent{\bf Acknowledgments.}
 The author was supported by
the National Natural Science Foundation of
China  (No. 11901430).

\end{document}